\magnification=\magstep1                          
\hsize=15truecm                                   
\hoffset=1.2cm                                    

\parskip=5truept plus 1truept minus 1truept       
\magnification \magstep1                          
\hsize=15truecm                                   
\hoffset=1.cm                                     

\parskip=5truept plus 1truept minus 1truept       
\baselineskip=20truept plus 1truept minus 1truept 

\font\grande=cmr10 scaled \magstep2               
\font\note=cmr8                                   
\font\maiuscoletto=cmcsc10 scaled \magstep1       

\def\titolo#1#2#3#4{{\centerline {\bf \grande {#1}}}\medskip
{\centerline {\bf \grande {#2}}}\medskip {\centerline {\bf \grande
{#3}}}\medskip {\centerline {\bf \grande {#4}}}} 
\def\abstract#1{\noindent{\bf Abstract.} {\narrower {\note #1}\par}}

\font\tenmsb=msbm10
\font\sevenmsb=msbm7
\font\fivemsb=msbm5
\newfam\msbfam
\textfont\msbfam=\tenmsb
\scriptfont\msbfam=\sevenmsb
\scriptscriptfont\msbfam=\fivemsb
\def\Bbb{\fam\msbfam }

\def\NN{{\Bbb N}}

\def\RR{{\Bbb R}}

\def\bull{{\vrule height.9ex width.8ex depth-.1ex}}       

\def\supp{\mathop{\rm supp} }

\newdimen\shift \newbox\leftbox 
\def\newitem#1{\par\setbox\leftbox=\hbox{#1}\shift=\wd\leftbox
\def\rightshift{\hskip -\parindent\hskip\shift\hangindent\shift}
\everypar={\rightshift}
\leavevmode\hskip -\shift\box\leftbox}
\def\endnewitem{\everypar={}}                     

\newdimen\rientrodescr \rientrodescr=20pt
 
\phantom{a} \vskip2truecm
\titolo {A MULTIPLICITY RESULT} {FOR THE LINEAR} {SCHR\"ODINGER-MAXWELL
EQUATIONS} {WITH NEGATIVE POTENTIAL}

\vskip2truecm                                      

\centerline {\maiuscoletto Giuseppe Maria Coclite }
\centerline  {\note S. I. S. S. A., via Beirut 2-4, Trieste 34014, Italy}
\centerline  {\note e-mail: coclite@sissa.it}

\vskip2truecm                                      

\abstract {In this paper it is proved the existence of a sequence of radial
solutions with negative energy of the linear Schr\"odinger-Maxwell equations
under the action of a negative potential.}

\vskip2truecm
\centerline {Ref. S.I.S.S.A. 20/2001/M (March 2001) }
\vskip5truecm

 \phantom{Introduzione}                                     
 \noindent {\bf 1. Introduction}
\bigskip

\bigskip
In this paper we study the interaction between the electomagnetic field
and the wave function related to a quantistic non-relativistic charged
particle, that is described by the Schr\"odinger equation.

In [2, 3, 11] it has been studied the case in which the electromagnetic field
is assigned. Here we shall assume  that the unknowns of the problem are both
the wave function $\> \psi = \psi (x,t) \>$ and the gauge potentials $\>
\varphi = \varphi (x,t)\>$ and $\> {\bf A}={\bf A}(x,t) \>$ related to the
electromagnetic fields $\>{\bf E}, \>{\bf H}\>$ by the equations $${\bf
E}=-{1\over c} {\partial{\bf A}\over {\partial t }}-  \nabla \varphi,\quad
{\bf H} =\nabla \times {\bf A}. $$

Such a situation has been studied by Benci and Fortunato (cfr. [5]) in the
case in which the charged particle "lives" in a space region $\> \Omega,\>$
which is bounded. Here we want to analyze the case in which $\> \Omega =
\RR^3.\>$ Moreover we assume that there is an external field deriving from a
potential $\> - V ( x ). \>$ We consider the electrostatic case, namely we
look for potentials $\> \varphi \>$ and $\> {\bf A} \>$ which do not depend on
time $\> t: \>$ $$ \varphi = \varphi (x),\quad {\bf A}={\bf A}(x),\qquad x \in
\RR^3,$$  and for standing wave function $$\psi (x,t) = u(x) e^{i \omega
t},\qquad x\in \RR^3 ,\> t\in \RR,$$ where $\>\omega \in \RR\>$ and $\>u\>$ is
real valued. In this situation we can assume $${\bf A}=0. $$ It can be shown
(cfr. [5]) that $\> \varphi, \>\omega \>$ and $\>u\>$ are related by the
equations  $$\cases{\displaystyle  -{1\over 2}\Delta u-\varphi u-V(x) u=\omega
u,\quad  & ${\rm in} \>\>\RR^3,$\cr  {}& {}\cr
 \Delta\varphi =4\pi u^2, & ${\rm in} \>\>\RR^3,$\cr } \eqno(1) $$
where $\> V : \RR^3 \rightarrow \RR\>$ is a radial positive map, which is the
potential of the external action. We shall assume:
\newitem {\hbox to .9truecm {\phantom{i}$ (V_1) $\hfill }} 
$\> V \>$ is continuous in $\> \RR^3 \, \backslash  \, \{ 0\}; \>$
\smallskip
\newitem {\hbox to .9truecm {\phantom{i}$ (V_2) $\hfill }}
 $\> V \in L^{3 \over 2} (\{ \vert x \vert \le 1 \});\>$
\smallskip \newitem {\hbox to .9truecm {\phantom{i}$ (V_3) $\hfill }}
 $\> \lim \limits_{\vert x \vert \to + \infty} V(x) =0;\>$
\smallskip
\newitem {\hbox to .9truecm {\phantom{i}$ (V_4) $\hfill }}
$\>\lim\limits_{\vert x \vert \rightarrow +\infty} x^2 V(x)= +\infty.\>$
\endnewitem
\bigskip
\noindent Observe that the coulumbian potential,
that is the most physically interesting one,
satisfies $\> (V_1),\> (V_2), \> (V_3)\>$  and $\> (V_4) \>$ (cfr. [13; 14]).

The equations in (1) have a variational structure infact they are the
Euler-Lagrange equations related to the functional: 

$$ F_{\omega}(u,\varphi)=$$
$$={1\over 4}\int_{\RR^3}\vert\nabla
u\vert^2dx- {1\over 2}\int_{\RR^3}\varphi u^2 dx -
{1\over 16 \pi}\int_{\RR^3}\vert\nabla \varphi \vert^2dx -{1\over 2}\int_{\RR^3} V(x) u^2 dx -
 {\omega\over 2}\int_{\RR^3}\vert u\vert^2dx, $$

\noindent this functional is strongly indefinite; this means that $\>
F_{\omega} \>$ is neither  bounded from below nor from above and this
indefinitness cannot be removed by a compact perturbation. Moreover $\>
F_{\omega} \>$ does not exibit symmetry properties. By a suitable variational
principle we are reduced to study an even functional which does not exibit the
same indefinitness of  $\> F_{\omega}. \>$ The main result of this paper is the
following. \bigskip

\bigskip
\noindent{\maiuscoletto Theorem 1 \phantom{a}} {\it Let} $\>V\>$ {\it
satisfy} $\> (V_1),\> (V_2), \> (V_3)\>$ {\it and} $\> (V_4) \>$ {\it then 
for all} $\> \omega <0 \>$ {\it problem} (1) {\it has 
infinitely many solutions} $\> \{(u_k,\varphi_k) \}_{k \in \NN}\> $  {\it with}
$\> u_k \in H^1 (\RR^3),\>$    
$$\int_{\RR^3} \vert
\nabla \varphi_k \vert^2 dx< \infty$$
{\it and such that} $\> F_{\omega} (u_k , \varphi_k) < - {\displaystyle {\omega
\over 2}}. \>$
\bigskip

\noindent The case in which $\> V\>$ is radially decreasing and belongs to $\>
L^p (\RR^3),\>$ with $\> {\displaystyle {3 \over 2}} <p< \infty, \>$ is
investigated in [9, Capitolo 6] and the nonlinear case is studied in [10].
Finally we recall that the Maxwell equations coupled with nonlinear
Klein-Gordon equation and with Dirac equation have been studied respectilvely
in [6; 12].
\bigskip

\bigskip
\noindent {\bf 2. The Variational Principle}
\bigskip

In this section we shall prove a variational principle which permits to reduce (1) to the study of the critical points of an even functional, which is not strongly indefinite. To this end we need some thecnical preliminaries.

We define the space $\> {\bf\cal D}^{1,2}(\RR^3) \> $ as the closure of 
 $\> C^{\infty }_0(\RR^3)\>$ with respect to the norm
 $$ \Vert u\Vert_{{\cal D}^{1,2}}=
\Big( \int_{\RR^3}\vert \nabla u\vert^2dx \Big)^ {1/2}.$$ 
 The following lemma holds (cfr. [7, Theorem 2.4]):
\bigskip

\bigskip
\noindent{\maiuscoletto Lemma 2 \phantom{a}}{\it For all} $\> \rho\in
L^1(\RR^3)\cap L^r(\RR^3),\>$ {\it with
 $\> {\displaystyle {6 \over 5} }<r\le 2,\> $
there exists only one }
 $\> \varphi \in {\cal D}^{1,2}(\RR^3)\>$ {\it such that}
 $\> \Delta\varphi =\rho.\>$
 {\it Moreover there results}
 $$ \Vert\varphi\Vert^2_{{\cal D}^{1,2}}\le c(\Vert\rho\Vert^2_{L^1}+
 \Vert\rho\Vert^2_{L^r} ) $$
 {\it and the map}
 $\> \rho\in L^1(\RR^3)\cap L^r(\RR^3)\mapsto\varphi=\Delta^{-1}(\rho)\in
 {\cal D}^{1,2}(\RR^3) \>$
 {\it is continuous.}
\bigskip

\bigskip
\noindent By Lemma 2 and Sobolev inequality, for any given $\>u \in H^1
(\RR^3)\>$ the second equation of (1) has the unique solution  $$ \varphi
=4\pi \Delta^{-1}u^2\>(\in {\cal D}^{1,2}(\RR^3)). $$ For this reason we can
reduce (1) to  
$$ -{1\over 2}\Delta u-4\pi(\Delta^{-1} u^2) u-V(x)u=\omega  u,\quad {\rm
in} \>\>\RR^3 \eqno(2)$$ 
Observe that (2) is the Euler-Lagrange equation
of the functional
 $$ J_{\omega}(u)={1\over 4}\int_{\RR^3}\vert\nabla u\vert^2dx+
 \pi\int_{\RR^3}\vert\nabla\Delta^{-1} u^2\vert^2dx-{1 \over 2} \int_{\RR^3}
V(x)u^2 dx -{\omega \over 2} \int_{\RR^3}
u^2 dx $$ Now we set 
$$ H^1_r(\RR^3):=\{u\in H^1(\RR^3)\big\arrowvert\> u(x)=u(\vert x\vert),
 \>\>\> x\in\RR^3 \}. $$
Since 
$${\displaystyle {d  \over {d \lambda}} \Big(
\int_{\RR^3}\vert\nabla\Delta^{-1} (u+ \lambda v) \vert^2dx}\Big)
\Big\arrowvert_{\lambda =0}= -2  \int_{\RR^3}(\Delta^{-1} u\vert v)dx$$ easly
the following lemma holds.
 \bigskip

\bigskip
\noindent{\maiuscoletto Lemma 3 \phantom{a}} {\it For all} $\> \omega \in \RR
\>${\it there results:} \smallskip
\newitem {\hbox to .7truecm {\phantom{ii}i)\hfill }} 
 $\> J_{\omega} \>$ {\it is even;}
\smallskip
\newitem {\hbox to .7truecm {\phantom{i}ii)\hfill }}
 $\> J_\omega \>$ {\it is} $\> C^1 \>$ {\it on} $\> H^1 (\RR^3) \>$ {\it and
its critical points are solutions of}
(2); 
\smallskip 
\newitem {\hbox to .7truecm {\phantom{i}iii)\hfill }}
 {\it any critical point of} $\> J_{\omega} \big\arrowvert_{H^1_r (\RR^3)} \>$
{\it is also a critical point of} $\> J_{\omega}. \>$
 \endnewitem
\bigskip

\bigskip
\noindent {\bf 3. Proof of Theorem 1}

\bigskip
We begin proving some lemmas.
\bigskip

\bigskip
\noindent{\maiuscoletto Lemma 4 \phantom{a}} {\it Let} $\>V \>$ {\it satisfy}
$\> (V_1),\> (V_2)\>$ {\it and } $\> (V_3)\>$ {\it then for all} $\> \omega <0
\>$ {\it the functional} $\> J_{\omega} \>$
{\it is weakly lower semicontinuous in} $\> H^1_r (\RR^3). \>$ {\it Precisely}
$$u \in H^1_r (\RR^3) \mapsto \int_{\RR^3}\vert \nabla u\vert^2dx - 2 \omega
\int_{\RR^3}u^2dx$$ {\it is weakly lower semicontinuous and } 
$$u \in H^1_r (\RR^3) \mapsto \int_{\RR^3}\vert\nabla\Delta^{-1}
u^2\vert^2dx,$$
 $$u \in H^1_r (\RR^3) \mapsto \int_{\RR^3}V(x) u^2dx$$ {\it are
weakly continuous.}  \bigskip

\noindent {\maiuscoletto Proof. \phantom{a}} Let $\> \omega <0. \>$ By a
well known argument the functional
$$u \in H^1_r (\RR^3) \mapsto \int_{\RR^3}\vert \nabla
u\vert^2dx - 2 \omega \int_{\RR^3}u^2dx$$
 is weakly lower semicontinuous.

Prove that the functional 
$$u \in H^1_r (\RR^3) \mapsto \int_{\RR^3}\vert\nabla\Delta^{-1} u^2\vert^2dx$$
is weakly continuous. We just have to observe that the operator
$$Q: \> u \in H^1_r (\RR^3) \mapsto u^2 \in L^{6 \over 5} (\RR^3) \cap L^2
(\RR^3)$$ is compact, infact by the compact embeddings of $\> H^1_r
(\RR^3 ) \>$ (cfr. [8, Theorem A.I'; 16]) the operator:
 $$H^1_r
(\RR^3)\hookrightarrow \hookrightarrow  L^{12 \over 5} (\RR^3) \cap L^4
(\RR^3){\buildrel Q\over \longrightarrow} L^{6 \over 5} (\RR^3) \cap L^2
(\RR^3)$$  
is compact and by Lemma 2 the following one
$$\Delta^{-1}: \>  L^{6 \over 5} (\RR^3) \cap L^2
(\RR^3) \longrightarrow {\cal D}^{1,2}(\RR^3)$$
 is continuous.

Prove that the functional
$$u \in H^1_r (\RR^3) \mapsto \int_{\RR^3}V(x) u^2dx$$
is weakly continuous. Let $\> \{u_k \} \subset H^1_r (\RR^3) \>$ and $\> u
\in H^1_r (\RR^3)\>$  such that
$$u_k \rightharpoonup u \quad {\rm weakly \>in}\>\> H^1_r (\RR^3).$$  
Since 
$$u_k \rightharpoonup u \quad {\rm weakly \>in}\>\> L^2 (\RR^3),$$
there exists $\> C>0 \>$ such that
$$\Vert u_k \Vert_{L^2} \le C, \quad\quad\quad \Vert u \Vert_{L^2} \le C.$$
By $\> (V_3) \>$ for all $\> \varepsilon >0 \>$ there exists $\> R>0 \>$ such
that 
$$\vert x \vert \le R \>\> \Longrightarrow 0\le V (x) < {\varepsilon
\over {C^2}}$$  
then
$$\int_{\{ \vert x \vert \ge R \}}V(x) u_k^2dx <\varepsilon , \quad\quad 
\int_{\{ \vert x \vert \ge R \}}V(x) u^2dx < \varepsilon. \eqno(3)$$
By the Sobolev inequality clearly
$$u_k^2 \rightharpoonup u^2 \quad {\rm weakly \>in}\>\>  L^3 (\RR^3),$$
and by $\> (V_1) \>$ and $\> (V_2) \>$ there results
$$\int_{\{ \vert x \vert \le R \}}V (x) u_k^2dx \rightarrow
\int_{\{ \vert x \vert \le R \}}V(x) u^2dx.$$
Then by the previous and (3) we can conclude
$$\int_{\RR^3}V(x) u_k^2dx \rightarrow
\int_{\RR^3}V(x) u^2dx.$$
So we are done.\bull 
\bigskip

\bigskip
\noindent{\maiuscoletto Remark 5 \phantom{a}} Observe that only for $\> 3 \le
n <6 \>$ we are able to prove that the functional 
$$u \in H^1_r (\RR^n) \mapsto \int_{\RR^n}\vert\nabla\Delta^{-1} u^2\vert^2dx$$
is weakly continuous by using the compact embedding results for radial
solutions (cfr. [8, Theorem A.I'; 16]) and Lemma 2.
\bigskip

\bigskip
\noindent{\maiuscoletto Lemma 6 \phantom{a}} {\it Let} $\>V \>$ {\it satisfy}
$\> (V_1),\> (V_2)\>$ {\it and } $\> (V_3)\>$ {\it then for all} $\> \omega <0
\>$ {\it the functional} $\> J_{\omega} \>$
{\it is coercive in} $\> H^1_r (\RR^3),\>$ {\it i. e. for all sequence} $\> \{
u_k \} \subset H^1_r (\RR^3) \>$ {\it such that} 
$\> \Vert u_k \Vert_{H^1} \rightarrow +\infty \>$ {\it there results}
$\> \lim\limits_k J_{\omega}(u_k)=+ \infty. \>$
\bigskip

\noindent{\maiuscoletto Proof. \phantom{a}} Let $\> \omega <0.
\>$ Denote
$$B' =\{ u \in H^1_r (\RR^3) \big\vert \Vert u \Vert_{H^1} =1 \}.$$
Let $\> \{u_k \} \subset H^1_r (\RR^3) \>$ such that 
$$\Vert u_k \Vert_{H^1} \rightarrow +\infty.$$
Denote
$$u_k = \lambda_k {\tilde u}_k$$
with $\> \lambda_k \in \RR \>$ and $\> {\tilde u}_k \in B' \>$. We have 
$$ J_{\omega} (u_k ) = a_k\lambda^2_k+b_k\lambda^4_k-c_k\lambda^2_k +
d_k\lambda^2_k$$ 
with 
$$a_k = {1 \over 4}\int_{\RR^3}\vert \nabla \tilde u_k\vert^2dx \in [0,
{1 \over 4}],\quad\quad  b_k = \pi \int_{\RR^3}\vert\nabla\Delta^{-1} \tilde
u_k^2\vert^2dx \ge 0,$$ 
 $$ c_k = {1 \over 2} \int_{\RR^3}V(x) {\tilde u}_k^2dx \ge 0,
\quad\quad d_k = -{\omega \over 2} \int_{\RR^3}{\tilde u}_k^2dx \in [0,
\, - {\omega \over 2}].$$
 Observe that by Sobolev inequality $\> (V_1),\> (V_2)\>$ and $\> (V_3)\>$
there results $$2 c_k =
\int_{\{ \vert x \vert \le 1 \}} V (x) {\tilde u}_k^2dx + \int_{\{ \vert x
\vert > 1 \}}  V(x) {\tilde u}_k^2dx \le$$  $$ \le \Vert V \Vert_{L^{3 \over
2} (\{ \vert x \vert \le 1 \})} \Vert \tilde u_k \Vert_{L^6 }^2 +
\sup\limits_{\vert x \vert \ge 1} V(x) \Vert \tilde u_k \Vert_{L^2 }^2 \le$$
$$\le \big( C \Vert V \Vert_{L^{3 \over 2} (\{ \vert x \vert \le 1 \})} +
\sup\limits_{\vert x \vert \ge 1} V(x) \big) \Vert \tilde u_k \Vert_{H^1}^2 =
\big( C \Vert V \Vert_{L^{3 \over 2} (\{ \vert x \vert \le 1 \})} +
\sup\limits_{\vert x \vert \ge 1} V(x) \big),$$ where $\> C > 0 \>$ is the
Sobolev embedding constant.  Since $\> u \in H^1_r (\RR^3) \mapsto
\int_{\RR^3}\vert\nabla\Delta^{-1} u^2\vert^2dx \>$ is weakly continuous and
$\> B' \>$ is bounded in $\> H^1_r (\RR^3 ) \>$ there exists $\> \alpha >0 \>$
such that $\> b_k \ge \alpha >0.\> $ 
Then we can conclude that
$$\lim\limits_k J_{\omega} (u_k)=+ \infty ,$$
and so we are done. \bull
\bigskip

\noindent By a well known argument by the two previous lemma the following
holds.

\bigskip
\noindent{\maiuscoletto Lemma 7 \phantom{a}} {\it Let} $\>V \>$ {\it satisfy}
$\> (V_1),\> (V_2)\>$ {\it and } $\> (V_3)\>$ {\it then for all} $\> \omega <
0 \>$ {\it the functional} $\> J_{\omega} \>$
{\it is bounded from below in} $\> H^1_r (\RR^3).\>$ 
\bigskip

\bigskip
\noindent{\maiuscoletto Lemma 8 \phantom{a}} {\it Let} $\>V \>$ {\it satisfy}
$\> (V_1),\>  (V_2) \>$ {\it and } $\> (V_3)\>$ {\it then for all} $\> \omega
<0 \>$ {\it the functional}$\> J_{\omega} \big\arrowvert_{H^1_r (\RR^3)} \>$
{\it stisfies the Palais-Smale condition, i.e. any sequence} $\quad \{ u_k \}
\subset H^1_r (\RR^3) \>$  {\it such that $\> \{ J_\omega (u_k) \}  \>$ is
bounded and  $\> J_\omega(u_k)\>\big\arrowvert_{H^1_r (\RR^3)} '
\rightarrow\>0 \>$ contains a coverging subsequence.} 
\bigskip

\noindent{\maiuscoletto Proof. \phantom{a}} Let $\> \omega <0 \>$ and $\> \{
u_k \} \subset H^1_r (\RR^3) \>$ such that $\> \{ J_\omega (u_k) \}  \>$ is
bounded and  $\> J_\omega(u_k)\>\big\arrowvert_{H^1_r (\RR^3)} '
\rightarrow\>0. \>$ First of all observe that, by (iii) of Lemma 3, there
results $$J_\omega\>\big\arrowvert_{H^1_r (\RR^3)} ' (u) =0 \Longleftrightarrow
J_\omega ' (u)=0,$$ then we can suppose $$J_\omega ' (u_k)  \rightarrow\>0.$$
By Lemma 6 $\> \{ u_k \} \>$ is bounded in $\> H^1_r (
\RR^3 ), \>$ passing to a subsequence there exists $\> u \in H^1_r (\RR^3 )\>$
such that
 $$ u_k\>\rightharpoonup\>u \>\>\quad \quad{\rm weakly \> in}\>
H^1_r(\RR^3). \eqno (4)$$
Clearly there results
$$J_{\omega}' (u) = 0. \eqno(5)$$ 
We prove that
$$ u_k\>\rightarrow\>u \>\>\quad \quad{\rm in}\> H^1_r(\RR^3).$$
By Lemma 4 and (4) there results
$$\int_{\RR^3}\vert \nabla u_k \vert^2dx -2 \omega \int_{\RR^3} u_k^2dx=$$
$$= 2\big\langle J_{\omega} '(u_k) ,\, u_k \big\rangle - 8 \pi
\int_{\RR^3}\vert\nabla\Delta^{-1} u_k^2\vert^2 dx +2\int_{\RR^3}V(x)
u_k^2dx \rightarrow $$
$$ \rightarrow - 8 \pi \int_{\RR^3}\vert\nabla\Delta^{-1} u^2\vert^2dx +
2\int_{\RR^3}V(x) u^2dx =$$ $$=\int_{\RR^3}\vert \nabla u \vert^2dx -2 \omega
\int_{\RR^3} u^2 dx - 2\big\langle J_{\omega} '(u) ,\, u \big\rangle.$$ 
By (5) and since $\> \omega < 0 \>$ the thesis is proved.\bull
 \bigskip

\bigskip
\noindent{\maiuscoletto Remark 9 \phantom{a}} Since for all $\> \omega <0 \>$
the functional $\> J_{\omega} \>$ is bounded from below and satisfies the
Palais-Smale condition there exists al least the critical level $\>
\inf J_{\omega}. \>$ The assumption $\> (V_4 ) \>$ helps us to obtain the
multiplicity of the same ones. \bigskip

\bigskip
\noindent{\maiuscoletto Lemma 10 \phantom{a}} {\it Let} $\>V \>$ {\it satisfy}
$\> (V_1),\> (V_2), \> (V_3)\>$ {\it and} $\> (V_4) \>$ {\it then for
all} $\> k \in \NN \backslash \{ 0 \}, \>$ {\it there exist a
subspace} $\> V_k \subset \>$ $\> H^1_r (\RR^3) \>$ {\it of dimension} $\> k
\>$ {\it and} $\> \nu > 0 \>$ {\it such that}  $$\int\limits_{\RR^3} \Big( {1
\over 2} \vert \nabla u \vert^2 - V(x) u^2 \Big) dx \le -\nu,$$ {\it for all}
$\> u \in V_k \cap B,\>$ {\it where} $$B =\Big\{ u \in H^1_r (\RR^3 )
\Big\arrowvert \int\limits_{\RR^3}  \vert  u \vert^2 dx =1 \Big\}.$$  \bigskip

\bigskip
\noindent{\maiuscoletto Proof. \phantom{a}} Let $\> u \>$ a smooth map with
compact support such that
$$\int\limits_{\RR^3}  \vert  u \vert^2 dx =1,\quad \supp (u)\subset B_2 (0)
\backslash B_1 (0),$$ 
where
$$B_{\rho} ( x ) =\big\{ y\in \RR^3 \big\arrowvert \vert x - y \vert < \rho
\big\},\quad\quad x \in \RR^3, \>\> \rho > 0.$$ Denote
$$u_{\lambda} (x)=\lambda^{3 \over 2} u(\lambda x), \quad\quad \lambda >0, \> x
\in \RR^3,$$ 
and
$$A_\lambda = B_{2 \over \lambda}(0) \backslash B_{1 \over \lambda}
(0),\quad\quad \lambda >0,$$ there results 
$$\int\limits_{\RR^3}  \vert  u \vert^2 dx =\int\limits_{\RR^3} 
\vert u_{\lambda}  \vert^2 dx =1, \quad \supp (u_{\lambda})\subset A_\lambda .$$ 
By $\> (V_1 ) \>$ we have $$\int\limits_{\RR^3} \Big( {1 \over 2}
\vert \nabla u_{\lambda} \vert^2 - V(x) u_{\lambda}^2 \Big)dx =
\int\limits_{\RR^3} \Big( \lambda^2 {1 \over 2} \vert \nabla u \vert^2 -
V \big( {x \over \lambda } \big) u^2 \Big)dx  \le $$ $$\le \lambda^2
\int\limits_{\RR^3} {1 \over 2} \vert \nabla u \vert^2 dx - 
\inf\limits_{{\supp u} \over \lambda} V \le \lambda^2 \int\limits_{\RR^3} {1
\over 2} \vert \nabla u \vert^2 dx -\inf\limits_{A_{\lambda}} V =$$
$$=\lambda^2 \int\limits_{\RR^3} {1 \over 2} \vert \nabla u \vert^2 dx - V(
x_\lambda ),$$ where $\> x_\lambda \>$ belongs to the clousure of $\> A_\lambda
\>$ and $\> V( x_\lambda ) =\displaystyle {\inf\limits_{A_{\lambda}}} V.\>$ By
$\> (V_3 ) \>$ and $\> (V_4 ) \>$ there exists $\> {\lambda}_0 >0 \>$ such
that  $$\int\limits_{\RR^3} \Big({1 \over 2} \vert \nabla u_{{\lambda}_0}
\vert^2 - V(x) u_{{\lambda}_0}^2 \Big)dx <0.$$ Let $\> k \in \NN \backslash \{
0 \} \>$ and $\> u_1, \, u_2, \, \ldots , \,  u_k \>$ smooth maps with compact
supports such that 
 $$\int\limits_{\RR^3}  \vert  u_i \vert^2 dx =1,\quad \supp
(u_i)\subset B_{2i} (0) \backslash B_i (0) , \quad i = 1, \, 2, \, \ldots
, \, k.$$ 
Using an analogous argument we are able to find $\> \lambda_1, \, \lambda_2,
\, \ldots , \,  \lambda_k >0\>$ such that
 $$\int\limits_{\RR^3} \Big( {1 \over
2} \vert \nabla u_{i_{{\lambda}_i}} \vert^2 - V(x) u_{i_{{\lambda}_i}}^2
\Big)dx <0, \quad i = 1, \, 2, \, \ldots , \, k.$$
Let
 $$0 < \bar \lambda < \min \{
\lambda_1, \, \lambda_2, \, \ldots , \,  \lambda_k \}$$ and $\> V_k
\>$ the subspace spanned by $\> u_{1_{\bar \lambda}}, \, u_{2_{\bar
\lambda}}, \, \ldots , \, u_{k_{\bar \lambda}}.\>$ 
Since the supports of this maps are pairwise disjoint
$\> V_k \>$ has dimension
$\> k. \>$ 
Since for all $\>i = 1, \, 2, \, \ldots , \, k \>$ and $\> \lambda \le
\lambda_i \>$there results
$$\int\limits_{\RR^3} \Big( {1 \over 2} \vert
\nabla u_{i_\lambda} \vert^2 - V(x) u_{i_\lambda}^2 \Big) < 0$$ 
and $\> V_k \cap B \>$ is compact, the thesis is proved.\bull 
\bigskip

\bigskip
\noindent{\maiuscoletto Lemma 11 \phantom{a}}{\it Let} $\>V \>$ {\it satisfy}
$\> (V_1),\> (V_2), \> (V_3)\>$ {\it and} $\> (V_4) \>$ {\it then for all} $\>
\omega < 0 \>$ {\it the functional} $\> J_{\omega} \>$ {\it has infinitely many
critical points} $\> \{ u_k \}_{k \in \NN} \subset H^1_r ( \RR^3 ) \>$ {\it
such that} $\> J_{\omega} ( u_k) < -{\displaystyle {\omega \over 2}}. \>$
\bigskip

\bigskip
\noindent{\maiuscoletto Proof. \phantom{a}} Let $\> \omega < 0 \>$ and denote
$$c^{\omega}_k = \inf \{ \sup J_{\omega}(A) \big\vert A \in {\cal A}, \> \gamma
(A) \ge k \}, \quad k \in \NN \backslash \{ 0 \},$$
with
 $${\cal A} = \{ A \subset H^1_r (\RR^3 ) \big\vert A \> {\rm closed ,\>
symmetric \> and \>} 0 \not\in A \}$$ 
and $\> \gamma \>$ is the Genus (cfr. e. g. [1, Definition 1.1]). We have to
prove that  $\> c^{\omega}_k < -{\displaystyle {\omega \over 2}}\>$ for all
$\> k \in \NN. \>$ 
Let $\> k \in \NN \backslash \{ 0 \} \>$, by the previous lemma there
exist $\> V_k \, \subset \, H^1_r (\RR^3)\>$ subspace of dimension
$\> k \>$   and $\> \nu > 0  \>$ such that for all $\> u
\in V_k \cap B \>$ there results
$$ \int\limits_{\RR^3} \Big( {1 \over 2} \vert \nabla u \vert^2 - V (x)u^2
\Big) dx \le -\nu.$$
Let $\> \lambda >0 \>$ and define
  $$h_\lambda : V_k \cap B \longrightarrow
H^1_r (\RR^3), \quad h_\lambda (u)= \lambda^{1 \over 2} u.$$ 
Fixed $\> u \in V_k \cap B\>$ and $\> 0 < \lambda < 1\>$ there results
  $$J_{\omega} ( h_\lambda (u)) \le  -{\lambda \over 2} \nu + c \lambda^2 -
{\omega \over 2} \lambda \le -{\lambda \over 2} \nu + c \lambda^2 -
{\omega \over 2} , \eqno(6)$$
 where $\> c
\>$ is a positive constant. Then there exists $\> 0 < \bar \lambda < 1 \>$ such
that for all  $\> u \in V_k \cap B \>$ there results $\> J_{\omega} (
h_{\bar\lambda} (u)) < -{\displaystyle {\omega \over 2}}.\>$ Since $\>
h_{\bar\lambda} \>$ is continuous, odd and $\> 0 \not\in V_k \cap B \>$ we
have $$ h_{\bar\lambda} (V_k \cap B) \in {\cal A}. \eqno(7)$$ Since $\>V_k
\cap B \>$ is compact, (6) and (7) we have
 $$\inf J_{\omega} \le c^{\omega}_k \le \sup J_{\omega}
(h_{\bar\lambda} (V_k \cap B)) <- {\omega \over 2}.$$ 
By Lemma 8 (cfr. [15, Theorem 9.1; 4]) there exists $\> \{ u_k \} \subset B\>$
sequence of critical points of $\> J_{\omega} \>$ such
that $\> J_{\omega}(u_k) =c^{\omega}_k < -{\displaystyle {\omega \over
2}}.\>$   So we are done.\bull 
\bigskip

\bigskip
\noindent{\maiuscoletto Proof of Theorem 1 \phantom{a}} Since 
$$F_{\omega} (u, 4 \pi \Delta^{-1} u^2 ) = J_{\omega} (u)$$
for all $\> \omega \in \RR \>$ and $\> u \in H^1 (\RR^3), \>$ by Lemma 3 and
the previous one the thesis is done.\bull
 \bigskip
\vskip10truecm

\centerline {\bf References }
\bigskip

\medskip
\newitem {\hbox to .8truecm {[1]} \hfill }{\maiuscoletto A. Ambrosetti,
 P. H. Rabinowitz,}
 {\it Dual variational methods in critical point theory and applications,}
 J. Funct. Analysis {\bf 14} (1973), 349-381. 
\endnewitem

\medskip
\newitem {\hbox to .8truecm {[2]} \hfill }{\maiuscoletto J. Avron,
 I. Herbsit, B. Simon,}
 {\it Schr\"odinger operators with magnetic fields I. General interaction,}
 Duke Math. J. {\bf 45} (1978), 847-883. 
\endnewitem

\medskip
\newitem {\hbox to .8truecm {[3]} \hfill }{\maiuscoletto J. Avron,
 I. Herbsit, B. Simon,}
 {\it Schr\"odinger operators with magnetic fields I. Atoms in homogeneous
 magnetic fields,} Comm. Math. Phys. {\bf 79} (1981), 529-572. 
\endnewitem

\medskip
\newitem {\hbox to .8truecm {[4]} \hfill }{\maiuscoletto P. Bartolo, V.
Benci, D. Fortunato, }
{\it Abstract critical point theorems and applications to some nonlinear
problems with "strong" risonance at infinity,} Nonlinear Analysis T.M.A.
 {\bf 7 } (1983) 981-1012.
\endnewitem

\medskip
\newitem {\hbox to .8truecm {[5]} \hfill }{\maiuscoletto V. Benci,
 D. Fortunato,}
 {\it An eigenvalue problem for the {\phantom \-} Schr\"odinger Maxwell
 equations,}
 Topological Methods in Nonlinear Analysis {\bf 11 } (1998), 283-293.
\endnewitem

\medskip
\newitem {\hbox to .8truecm {[6]} \hfill }{\maiuscoletto V. Benci, D.
Fortunato,}
 {\it The nonlinear Klein-Gordon equation coupled with the Maxwell equations, }
  to appear on Nonlinear Analysis T. M. A., Proceedings W. C. N. A. 2000.
\endnewitem

\medskip
\newitem {\hbox to .8truecm {[7]} \hfill }{\maiuscoletto V. Benci, D.
Fortunato, A. Masiello, L. Pisani,  }
 {\it Solitons and electromagnetic field, }
  Math. Z. {\bf 232} (1999), 73-102.
\endnewitem
\medskip
\newitem {\hbox to .8truecm {[8]} \hfill }{\maiuscoletto H. Beresticki, P. L.
Lions, }  {\it Nonlinear scalar field equations. I Existence of a ground
state, }  Arch. Rational Mech. Anal. {\bf 82} (1983), 313-346.
\endnewitem

\medskip
\newitem {\hbox to .8truecm {[9]} \hfill }{\maiuscoletto G. M. Coclite,}
 {\it Metodi Variazionali Applicati allo Studio delle Equazioni di
Schr\"odinger-Maxwell,} thesis, University of Bari (1999). \endnewitem

\medskip
\newitem {\hbox to .8truecm {[10]} \hfill }{\maiuscoletto G. M. Coclite,}
 {\it A Multiplicity Result for the  Nonlinear Schr\"odinger-Maxwell
Equations,} to appear on Communications on Applied Analysis.
\endnewitem

\medskip
\newitem {\hbox to .8truecm {[11]} \hfill }{\maiuscoletto J. M. Combes, 
 R. Schrader, R. Seiler,}
 {\it Classical bounds and limits for energy distributions of Hamiltonian
 operators in electromagnetic fields,} Ann. Phys. {\bf 111} (1978), 1-18.
\endnewitem

\medskip
\newitem {\hbox to .8truecm {[12]} \hfill }{\maiuscoletto M. J. Esteban, 
 V. Georgiev, E. Sere,}
 {\it  Stationary solutions of the Maxwell-Dirac and the Klein-Gordon-Dirac
equations,} Calc. Var. {\bf 4} (1996), 265-281. 
\endnewitem

\medskip
\newitem {\hbox to .8truecm {[13]} \hfill }{\maiuscoletto P. L. Lions,  }
 {\it Solutions of Hartree-Fock equations for Coulumb systems, }
 Commun. Math. Phys. {\bf 109} (1987), 33-97.
\endnewitem

\medskip
\newitem {\hbox to .8truecm {[14]} \hfill }{\maiuscoletto P. L. Lions, }
 {\it Some remarks on Hartree equations, }
 Nonlinear Analysis T. M. A. {\bf 5} (1981), 1245-1256.
\endnewitem

\medskip
\newitem {\hbox to .8truecm {[15]} \hfill }{\maiuscoletto P. H. Rabinowitz, }
 {\it Minimax methods in critical point theory with apllications to
differential equations, } AMS (1986).
\endnewitem

\medskip
\newitem {\hbox to .8truecm {[16]} \hfill }{\maiuscoletto W. A. Strauss,  }
 {\it Existence of solitary waves in higher dimensions, }
 Commun. Math. Phys. {\bf 55} (1977), 149-162.
\endnewitem

\end